\documentclass[12pt]{article}  
\usepackage{amssymb}

\newtheorem{corollary}{Corollary}
\newtheorem{proposition}{Proposition}

\newenvironment{definition}
{\smallskip\noindent{\bf Definition\/}:}{\smallskip\par}

\newenvironment{remarks}
{\smallskip\noindent{\bf Remarks\/}.}{\smallskip\par}
\newenvironment{proof}
{\noindent{\em Proof\/}.}{{ $\Box$}\smallskip\par}

\newcommand{\CC}{{\Bbb C}}

\newcommand{\PP}{{\Bbb P}}

\newcommand{\ZZ}{{\Bbb Z}}

\title{A filtration defined by arcs on a variety}
\author{W.Ebeling and S.M.Gusein-Zade
\thanks{Partially supported by the DFG-programme ''Global methods in
complex geometry'' (Eb 102/4--2), grants RFBR--01--01--00739,
INTAS--00-259.
}
}
\date{}

\begin{document}

\maketitle

\begin{abstract} We define a natural filtration on the ring ${\cal O}_{V,0}$ of germs
of functions on a germ of a complex analytic variety $(V,0)$ related with the
geometry of arcs on the variety and describe some properties of it.
\end{abstract}

Let $(V,0)$ be a germ of a complex analytic variety and let ${\cal O}_{V,0}$ be the ring of germs
of functions on it. There is a natural filtration on the ring ${\cal O}_{V,0}$ defined by powers
of the maximal ideal. We define another natural filtration related with the geometry of arcs on
$(V,0)$ and describe some properties of it. 

An arc $\phi$ on $(V,0)$ is a germ of a complex analytic mapping $\phi: (\CC,0) \to (V,0)$. For
a function $g \in {\cal O}_{V,0}$, its order $v_\phi(g)$ on the arc $\phi$ is defined as the order
of the composition $g \circ \phi$, i.e., as the power of the first non-vanishing term in the power
series expansion $g \circ \phi (\tau) = a \tau ^{v_\phi(g)} + \ldots $, $a \neq 0$ (if $g \circ
\phi \equiv 0$, then $v_\phi(g)$ is assumed to be equal to $+\infty$). Let $\widehat{v}(g):=
\min_{\{ \phi \}} v_\phi(g)$. One has $\widehat{v}(g) < \infty$ for $g \neq 0$. One can easily
see that $\widehat{v} : {\cal O}_{V,0}\setminus \{0\} \to \ZZ_{\geq 0}$ is a valuation.

\begin{definition} The {\em arc filtration} 
${\cal O}_{V,0}=F_0 \supset F_1 \supset F_2 \supset \ldots$ on
the ring ${\cal O}_{V,0}$ is the filtration by the ideals $F_i= \{ g \in {\cal O}_{V,0} \,
|\, \widehat{v}(g) \geq i\}$. 
\end{definition}

On a germ of a smooth variety the arc filtration coincides with the filtration by powers
of the maximal ideal.

Suppose that $(V,0)$ has an isolated singularity and let $\pi : (X, D) \to (V,0)$ be a resolution
of it where $D=\bigcup_{i=1}^r E_i$ is a normal crossing divisor and $E_i$ ($i=1, \ldots , r$)
are its irreducible components.  Let $v_i(g)$ ($g \in {\cal O}_{V,0}$) be the order of the
lifting $\widetilde{g}= g \circ \pi$ of the function $g$ to the space $X$ of the resolution along
the component $E_i$ ($i=1, \ldots , r$). This defines a collection of valuations $\{ v_1, \ldots
, v_r \}$ on the ring ${\cal O}_{V,0}$ and the corresponding multi-indexed filtration on it:
$J(\underline{v})=\{ g \in {\cal O}_{V,0} \, | \, \underline{v}(g) \geq \underline{v} \}$,
$\underline{v}=(v_1, \ldots , v_r)$, $\underline{v}(g)=(v_1(g), \ldots , v_r(g))$ (see, e.g.,
\cite{CDG3}).

\begin{proposition} \label{Prop1}
$$\widehat{v}(g)= \min_{1Ê\leq i \leq r} v_i(g).$$
\end{proposition}

\begin{proof} This follows from the fact that the set of arcs on $(V,0)$ coincides with
the set of arcs on $(X,D)$ and the order of the function $g$ along the component $E_i$ of the
exceptional divisor $D$ is equal to its order along a generic smooth curve transversal to $E_i$.
\end{proof}

Let 
$$P_{V,0}(t)= \sum_{i=0}^\infty \dim (F_i/F_{i+1}) \cdot t^i$$
be the Poincar\'e series of the arc filtration on the ring ${\cal O}_{V,0}$.

\begin{sloppypar}

In \cite{CDK, CDG3} there was defined a (generalized) Poincar\'e series 
$P_{\{v_i\}}(t_1, \ldots , t_r)$ corresponding to the collection $\{v_i\}$ of valuations.
Strictly speaking there it is defined for the ring of functions on a curve or on a surface,
however, there is no difference to the general case.

\end{sloppypar}

For a formal power series $P(t_1, \ldots , t_r) \in \ZZ[[\, \underline{t} \, ]]$ 
($\underline{t}=(t_1,
\ldots , t_r)$), let its {\em reduction} $\overline{P}(t)$ be the formal power series in one
variable
$t$ obtained from the series $P(t_1, \ldots , t_r)$ by substituting each monomial
$\underline{t}^{\underline{v}}=t_1^{v_1} \cdots t_r^{v_r}$ in it by the monomial $t^{\min v_i}$.
For example, if $P(t_1,t_2)=1+t_1t_2^2+t_1^2t_2$, then $\overline{P}(t)=1+2t$.

\begin{proposition} \label{Prop2}
Let the germ $(V,0)$ be irreducible. Then the Poincar\'e series $P_{V,0}(t)$ of the arc
filtration coincides with the reduction $\overline{P_{\{v_i\}}}(t)$ of the Poincar\'e series  
$P_{\{v_i\}}(t_1, \ldots , t_r)$ of the collection $\{ v_i\}$ of the divisorial valuations.
\end{proposition}

\begin{proof}
A convenient way to see this is to express both Poincar\'e series as certain integrals with
respect to the Euler characteristic over the projectivization $\PP{\cal O}_{V,0}$ of the space
${\cal O}_{V,0}$ of germs of functions on $(V,0)$. The definition of such an integral 
can be found, e.g., in
\cite{CDG2,CDG3}. It is inspired by
the notion of motivic integration (see, e.g., \cite{DL}) and in some sense dual to it.  Let
$\ZZ[[\, \underline{t} \, ]]$ 
be the abelian group (with respect to addition) of formal power series in $t_1,
\ldots, t_r$. A monomial $\underline{t}^{\underline{v}(g)}$ can be considered as a
function on the projectivization $\PP{\cal O}_{V,0}$ of the space ${\cal O}_{V,0}$ with values in
the group $\ZZ[[\, \underline{t} \, ]]$. From
\cite{CDG3} it follows that
\begin{equation}
P_{\{v_i\}}(t_1, \ldots , t_r)= \int_{\PP{\cal O}_{V,0}} 
\underline{t}^{\underline{v}(g)} d\chi. 
\end{equation}
The proof in \cite{CDG3} is formulated for surface singularities, but it can be easily extended
to an arbitrary collection of finitely determined valuations $\{ v_i\}$. (A valuation $v$
is {\em finitely determined} if for any $i \geq 0$ there exists $N \geq 0$ such that $\{g \in
{\cal O}_{V,0} \, | \, v(g) \geq i\} \supset \frak{m}^N$.) For the arc valuation $\widehat{v}$ it
gives
\begin{equation}
P_{V,0}(t)= \int_{\PP{\cal O}_{V,0}} 
t^{\widehat{v}(g)} d\chi. 
\end{equation}
The formulae (1) and (2) and Proposition~\ref{Prop1} imply Proposition~\ref{Prop2}.
\end{proof}

The statement of Proposition~\ref{Prop3} does not hold in general for reducible germs. 

\begin{corollary} The reduction $\overline{P_{\{v_i\}}}(t)$ of the Poincar\'e series  
$P_{\{v_i\}}(t_1, \ldots , t_r)$ of the collection $\{ v_i\}$ of the divisorial valuations
defined by a resolution does not depend on the resolution. $\Box$
\end{corollary}

For rational surface singularities, the Poincar\'e series $P_{\{v_i\}}(t_1, \ldots , t_r)$ (for
the minimal resolution) were computed in \cite{CDG3}. Using Proposition~\ref{Prop2}, these
results yield the following statement.

\begin{proposition} \label{Prop3}
For the rational double points one has the following list of the Poincar\'e series 
$P_{V,0}(t)$ of the arc filtration:
\begin{eqnarray*}
A_k & : & \frac{1-t^2}{(1-t)^3} \\
D_k & : & \frac{1-t^{k-1}}{(1-t)^2(1-t^{k-2})} \\
E_6 & : & \frac{1-t^4}{(1-t)(1-t^2)^2} \\
E_7 & : & \frac{1-t^6}{(1-t)(1-t^2)(1-t^4)} \\
E_8 & : & \frac{1-t^6}{(1-t^2)^2(1-t^3)}
\end{eqnarray*}
\end{proposition}

\begin{remarks} {\bf 1.} The computations, though formal and routine, are somewhat tedious. In
each step one has in general expressions which are not products of powers of cyclotomic
polynomials. Moreover, for some singularities (e.g., for $D_5$) they even start from one which is
not of this sort (see \cite{CDG3}). However, finally they lead to series which are such 
products.

\noindent {\bf 2.} For an arbitrary singularity $(V,0)$ the Poincar\'e series $P_{V,0}(t)$ of the
arc filtration must not be a product of powers of cyclotomic polynomials. For example, for the
monomial space curve given by $x \mapsto (x^3,x^4,x^5)$ one has
$P_{V,0}(t)=(1-t+t^3)/(1-t)$. However, for an irreducible {\em plane} curve singularity the
Poincar\'e series $P_{V,0}(t)$ coincides with the Poincar\'e series considered, e.g., in
\cite{CDG1} and therefore is of this type.

\noindent {\bf 3.} Let $(V,0)$ be an irreducible singularity with a good $\CC^\ast$-action (i.e.,
$0$ is in the closure of every orbit). Then, if the $\CC^\ast$-action is free outside of the
origin, the arc filtration coincides with the one defined by the natural grading on the ring
${\cal O}_{V,0}$ corresponding to the $\CC^\ast$-action. However, if the action is not free, this
is not the case.

\noindent {\bf 4.} One can see that for all rational double points the Poincar\'e series
$P_{V,0}(t)$ of the arc filtration (as a rational function of $t$) has degree $-1$ and a pole of
order $2$ at $t=1$. 
Moreover, all these series coincide with the Poincar\'e series with respect to the
quasihomogeneous grading of certain rational double points. Namely, this correspondence is the
following one:
$$A_k \mapsto A_1, \quad D_k \mapsto A_{k-2}, \quad, E_6 \mapsto A_3, \quad E_7 \mapsto A_5,
\quad E_8 \mapsto D_4.$$
We don't understand the meaning of this correspondence.
\end{remarks}

The fact that all singularities $A_k$ have one and the same Poincar\'e series
$P_{V,0}(t)$ is a particular case of the following statement (which can be easily proved).

\begin{proposition}
Let $f$ be a germ of a function in $n$ variables $z_1, \ldots , z_n$ which belongs to the square
$\frak{m}^2$ of the maximal ideal $\frak{m} \subset {\cal O}_{\CC^n,0}$ and let $(V,0)$ be the
hypersurface singularity in $(\CC^{n+2},0)$ $($with the coordinates $(z_1, \ldots, z_n,x,y)$$)$
defined by the equation
$f(z_1, \ldots, z_n) +x^2+y^2 =0$ $($the double suspension of the hypersurface singularity
$\{f=0\}$$)$. Then 
$$P_{V,0}(t)= \frac{1-t^2}{(1-t)^{n+2}}.$$
\end{proposition}

The main property which guarantees that $P_{V,0}(t)$ does not depend on $f$ in $\frak{m}^2$ for
the hypersurface $V=\{ f(z_1, \ldots, z_n) +x^2+y^2 =0 \}$ is that there are many smooth curves on
$V$.

\bigskip
\noindent Universit\"{a}t Hannover, Institut f\"{u}r Mathematik \\
Postfach 6009, D-30060 Hannover, Germany \\
E-mail: ebeling@math.uni-hannover.de\\

\medskip
\noindent Moscow State University, Faculty of Mechanics and Mathematics\\
Moscow, 119992, Russia\\
E-mail: sabir@mccme.ru

\end{document}